\documentclass[12pt,leqno,draft]{article}

\def\dist{\mathop{\rm dist}}
\def\osc{\mathop{\rm osc}}

\begin{document}

\title{Some remarks about metric spaces, 2}

\author{Stephen Semmes\thanks{This survey has been prepared in
connection with the workshop on discrete metric spaces and their
applications at Princeton, August, 2003. }}

\date{}

\maketitle

	Let $(M, d(x,y))$ be a metric space.  Thus $M$ is a nonempty
set, and $d(x, y)$ is a real-valued function defined for $x, y \in M$
such that $d(x, y) \ge 0$ and $d(x, y) = d(y, x)$ for all $x, y \in
M$, $d(x, y) = 0$ if and only if $x = y$, and
\begin{equation}
	d(x, z) \le d(x, y) + d(y, z)
\end{equation}
for all $x, y, z \in M$.  This last property is called the
\emph{triangle inequality}.

	Suppose that $f(x)$ is a real-valued function on $M$.  If $C$
is a nonnegative real number, then we say that $f(x)$ is
\emph{$C$-Lipschitz} if
\begin{equation}
	|f(x) - f(y)| \le C \, d(x, y)
\end{equation}
for all $x, y \in M$, which is equivalent to saying that
\begin{equation}
	f(x) \le f(y) + C \, d(x, y)
\end{equation}
for all $x, y \in M$.  Notice that a function is $0$-Lipschitz
if and only if it is constant.

	For instance, for each $p \in M$, the function $f_p(x) = d(x,
p)$ is $1$-Lipschitz.  More generally, if $A$ is a nonempty subset of
$M$, then the distance of a point $x$ in $M$ to $A$ is denoted
$\dist(x, A)$ and defined by
\begin{equation}
	\dist(x, A) = \inf \{d(x, a) : a \in A\},
\end{equation}
and one can check that this function is $1$-Lipschitz.  If $f_1$,
$f_2$ are two real-valued functions on $M$ which are $C_1$,
$C_2$-Lipschitz, respectively, and if $\alpha_1$, $\alpha_2$ are real
numbers, then $\max(f_1, f_2)$, $\min(f_1, f_2)$ are $C$-Lipschitz
with $C = \max (C_1, C_2)$, and $\alpha_1 \, f_1 + \alpha_2 \, f_2$ is
$C$-Lipschitz with $C = |\alpha_1| \, C_1 + |\alpha_2| \, C_2$.

	Now suppose that $C$ is a nonnegative real number and that $s$
is a positive real number.  A real-valued function $f(x)$ on $M$ 
is said to be \emph{$C$-Lipschitz of order $s$} if
\begin{equation}
	|f(x) - f(y)| \le C \, d(x, y)^s
\end{equation}
for all $x, y \in M$, which is again equivalent to
\begin{equation}
	f(x) \le f(y) + C \, d(x, y)^s
\end{equation}
for all $x, y \in M$.  As before, $f(x)$ is $0$-Lipschitz of order $s$
if and only if $f(x)$ is constant on $M$.

	When $0 < s < 1$, one can check that $d(x, y)^s$ is also a
metric on $M$, which defines the same topology on $M$ in fact.  The
main point in this regard is that the triangle inequality continues to
hold, which follows from the observation that
\begin{equation}
	(\alpha + \beta)^s \le \alpha^s + \beta^s
\end{equation}
for all nonnegative real numbers $\alpha$, $\beta$.  A real-valued
function $f(x)$ on $M$ is $C$-Lipschitz of order $s$ with respect to
the metric $d(x, y)$ if and only if $f(x)$ is $C$-Lipschitz of order
$1$ with respect to $d(x, y)^s$, and as a result when $0 < s < 1$ one
has the same statements for Lipschitz functions of order $s$ as for
ordinary Lipschitz functions.

	When $s > 1$ the triangle inequality for $d(x, y)^s$ does not
work in general, but we do have that
\begin{equation}
	d(x, z)^s \le 2^{s-1} \, (d(x, y)^s + d(y, z)^s)
\end{equation}
for all $x, y, z \in M$, because
\begin{equation}
	(\alpha + \beta)^s \le 2^{s-1} \, (\alpha^s + \beta^s)
\end{equation}
for all nonnegative real numbers $\alpha$, $\beta$.  Some of the usual
properties of Lipschitz functions carry over to Lipschitz functions of
order $s$, perhaps with appropriate modification, but for instance it
may be that the only Lipschitz functions of order $s$ when $s > 1$ are
constant.

	Of course Lipschitz functions of any order are continuous.
The Lipschitz conditions provide concrete quantitative versions of the
notion of continuity.  Let us point out that in general the product of
two functions which are Lipschitz of order $s$ may not be Lipschitz of
order $s$, but that this is the case if at least one of the functions
is bounded.

	In harmonic analysis one considers a variety of classes of
functions with different kinds of restrictions on size, oscillations,
regularity, and so on, and these Lipschitz classes are fundamental
examples.  In particular, it can be quite useful to have the parameter
$s$ available to adjust to the given circumstances.  There are also
other ways of introducing parameters to get interesting classes of
functions and measurements of their behavior.

	If $M$ is the usual $n$-dimensional Euclidean space ${\bf
R}^n$, with its standard metric, then one has the extra structure of
translations, rotations, and dilations.  If $f(x)$ is a real-valued
function on ${\bf R}^n$ which is $C$-Lipschitz of order $s$, $f(x -
u)$ is also $C$-Lipschitz of order $s$ for each $u \in {\bf R}^n$,
$f(\Theta(x))$ is $C$-Lipschitz of order $s$ for each rotation
$\Theta$ on ${\bf R}^n$, and $f(r^{-1} x)$ is $(C \, r^s)$-Lipschitz
of order $s$ for each $r > 0$.  In effect, one general metric spaces
we can consider classes of functions and measurements of their
behavior which have analogous features, even if there are not exactly
translations, rotations, and dilations.

	On Euclidean spaces there is the classical \emph{Fourier
transform}, and for instance smoothness of a function can be related
to the size of the Fourier transform in various ways.  With the
Fourier transform there are very precise versions of information at
different wavelengths, including very specific ranges of wavelengths.
As in the Heisenberg uncertainty principle, however, there is a
balance between details of location and details of ranges of
wavelengths.

	With simple measurements like $t^{-s} \osc(x, t)$, one has
some information about location and range of wavelengths, but not too
precisely for either one.  Quantities like these also make sense in
general settings, without a lot of fine structure as for Euclidean
spaces.  At the same time, one gets at information and structure which
is interesting in the classical case of Euclidean spaces as well as 
other situations.

	A basic notion is to consider various scales and locations
somewhat independently.  In this regard, if $f(x)$ is a real-valued
function on $M$, $x$ is an element of $M$, and $t$ is a positive real
number, put
\begin{equation}
	\osc(x, t) = \sup \{|f(y) - f(x)| : y \in M, d(y, x) \le t\}.
\end{equation}
We implicitly assume here that $f(y)$ remains bounded on bounded
subsets of $M$, so that this quantity is finite.

	Thus $f$ is $C$-Lipschitz of order $s$ if and only if
\begin{equation}
	t^{-s} \osc(x, t) \le C
\end{equation}
for all $x \in M$ and $t > 0$.  Now, instead of considering basic
Lipschitz conditions like these, one can also look at other kinds of
bounds for $t^{-s} \osc(x, t)$.  Moreover, one can consider other
kinds of local measurements of size and oscillation.

	Let us pause a moment and notice that
\begin{equation}
	\osc(w, r) \le \osc(x, t)
\end{equation}
when $d(w, x) + r \le t$.  Thus, 
\begin{equation}
	r^{-s} \, \osc(w, r) \le 2^s \, t^{-s} \, \osc(x, t)
\end{equation}
when $d(x, w) + r \le t$ and $r \ge t/2$.

	This is a kind of ``robustness'' property of these
measurements of local oscillation of a function $f$ on $M$.  In
particular, to sample the behavior of $f$ at essentially all locations
and scales, it is practically enough to look at a reasonably-nice and
discrete family of locations and scales.  For instance, one might
restrict one's attention to radii $t$ which are integer powers of $2$,
and for a specific choice of $t$ use a collection of points in $M$
which cover suitably the various locations at that scale.

	Instead of simply taking a supremum of some measurements of
local oscillation like this, one can consider various sums of discrete
samples of this sort.  This leads to a number of classes of functions
and measurements of their behavior.  One can adjust this further
by taking into account the relation of some location and scale
to some kind of boundaries, or singularities, or concentrations,
and so on.

	Of course one might also use some kind of measurement of sizes
of subsets of $M$.  This could entail volumes, or sizes in terms of
covering conditions, or measurements of capacity.  One can then
look at integrals of $f$ and its powers, integrals involving the
local oscillation numbers $\osc(x, t)$, sizes of sets where some
other measurements are large, etc.

	There are also many kinds of local measurements of oscillation
or size that one can consider.  As an extension of just taking
suprema, one can take various local averages, or averages of powers of
other quantities.  Of course one can still bring in powers of the
radius as before.

	Even if one starts with measurements of localized behavior
which are not so robust in the manner described before, one can
transform them into more robust versions by taking localized suprema
or averages or whatever afterwards.  Frequently the kind of overall
aggregations employed have this kind of robustness included in effect,
and one can make some sort of rearrangement to put this in starker
relief.  Let us also note that one often has local measurements which
can be quite different on their own, but in some overall aggregation
lead to equivalent classes of functions and similar measurements of
their behavior.

	There are various moments, differences, and higher-order
oscillations that can be interesting.  As a basic version of this, one
can consider oscillations of $f(x)$ in terms of deviations from
something like a polynomial of fixed positive degree, rather than
simply oscillations from being constant, as with $\osc(x, t)$.  This
can be measured in a number of ways.

	However, for these kinds of higher-order oscillations,
additional structure of the metric space is relevant.  On Euclidean
spaces, or subsets of Euclidean spaces, one can use ordinary
polynomials, for instance.  This carries over to the much-studied
setting of nilpotent Lie groups equipped with a family of dilations,
where one has polynomials as in the Euclidean case, with the
degrees of the polynomials defined in a different way using the
dilations.

	These themes are closely related to having some kind of
derivatives around.  Just as there are various ways to measure the
size of a function, one can get various measurements of oscillations
looking at measurements of sizes of derivatives.  It can also be
interesting to have scales involved in a more active manner, and in
any case there are numerous versions of ideas along these lines
that one can consider.

\end{document}